\documentclass[conference]{IEEEtran}

\usepackage{dcolumn}    
\usepackage{eurosym}
\usepackage{graphicx}
\usepackage{subfig}
\usepackage{amsmath}
\usepackage{multirow}
\usepackage{amsmath}
\usepackage{amsfonts}
\usepackage{amssymb}

\usepackage[normalem]{ulem}

\graphicspath{{Images/}}
\newcolumntype{d}[1]{D{.}{.}{#1}}
\usepackage{siunitx} 
\usepackage{fancyhdr,lipsum}

\fancypagestyle{kbmcsc}{%
   \fancyhf{} 
   
   \fancyhead[L]{SUBMITTED TO IEEE ISGT NA 2019}
}%

\makeatletter
\let\ps@IEEEtitlepagestyle\ps@kbmcsc
\makeatother

\ifCLASSINFOpdf
\else
\fi
\begin{document}
%
\title{The Role of Demand-Side Flexibility in Hedging Electricity Price Volatility in Distribution Grids}

\author{\IEEEauthorblockN{Shantanu Chakraborty, \\Remco Verzijlbergh,  Zofia Lukszo}
\IEEEauthorblockA{Faculty of Technology, Policy, and Management\\
Delft University of Technology\\
Delft, Netherlands 2628BX\\
Email: S.T.Chakraborty@tudelft.nl}
\and
\IEEEauthorblockN{Milos Cvetkovic}
\IEEEauthorblockA{Faculty of Electrical Engineering, \\Mathematics and Computer Science\\
Delft University of Technology\\
Delft, Netherlands 2628CD\\
Email: M.Cvetkovic@tudelft.nl}
\and
\IEEEauthorblockN{Kyri Baker}
\IEEEauthorblockA{Civil, Environmental and \\Architectural Engineering\\
University of Colorado at Boulder\\
Boulder, CO 80309\\
Email: Kyri.Baker@colorado.edu}}


%


\maketitle

\begin{abstract}
Locational Marginal Price (LMP) is a dual variable associated with supply-demand matching and represents the cost of delivering power to a particular location if the load at that location increases. In recent times it become more volatile due to increased integration of renewables that are intermittent. The issue of price volatility is further heightened during periods of grid congestion. Motivated by these problems, we propose a market design where, by constraining dual variables, we determine the amount of demand-side flexibility required to limit the rise of LMP. Through our proposed approach a price requesting load can specify its maximum willingness to pay for electricity and through demand-side flexibility hedge against price volatility. For achieving this, an organizational structure for flexibility management is proposed that exhibits the coordination required between the Distribution System Operator (DSO), an aggregator and the price requesting load. To demonstrate the viability of our proposed formulation, we run an illustrative simulation under infinite and finite line capacities. 
\end{abstract}


%
\IEEEpeerreviewmaketitle

\section{Introduction}\label{sec1}
Large scale integration of distributed energy resources, such as wind and solar, poses both advantages and challenges for power systems. While it positively contributes to the reduction of carbon footprint of electricity, one of the issues related to the variability of distributed energy resources is that of price volatility. Demand side uncertainty in addition to intermittency in generation supply can cause the prices to be volatile in electricity markets \cite{Astaneh2013}. Price spikes, which are an extreme of price volatility, occur when there are no price-responsive loads, enabling generators to exercise market power; or when there is limited supply in energy markets to satisfy peak loads or during periods of grid congestion \cite{Asadinejad2017}.\par

Price volatility due to intermittent nature of renewables is observed in different regions around the world such as in Spain \cite{Pereira2017}, Australia \cite{Higgs2015}, Denmark \cite{Rintamaki2017}, Germany \cite{Ketterer2014}, and Texas \cite{Woo2011}. In some locations supply-side flexibility options in the form of hydropower \cite{Pereira2017,Rintamaki2017} or cross-border power flow \cite{Rintamaki2017} is able to address price volatility caused by the increased integration of renewables. However, supply-side flexibility options are subject to grid congestions which could possibly lead to price spikes. In \cite{Lund2015} the adoption of demand-side flexibility in addition to supply-side is recommended to address issues such as market power of generators, reducing price spikes and price volatility, which in turn can support higher integration of renewables. In this paper, we are concerned with the demand-side flexibility, which can be availed from energy storage and demand response, and is defined as the ability of a power system to modify consumption in reaction to an external signal such as price. \par

To address the issue of price volatility and price spikes both regulatory and technology based approaches have previously been investigated. From a regulatory perspective research in \cite{Asadinejad2017} proposes an optimal combination of the time-of-use and incentive based demand response with a focus on customer comfort. Similarly, \cite{Mohajeryami2017} focuses on analyzing the consumer reaction and their willingness to shift demand under demand response programs. In contrast, the optimal size of energy storage as a demand-side flexibility resource for addressing price volatility in a renewable generation-dominant electricity market is investigated in \cite{Masoumzadeh2018}. \par

Hedging has also been used as an approach for addressing price volatility and price spikes. In \cite{Baatz2018} the research focuses on the potential of Energy Efficiency as a hedging mechanism, which is similar to long-term energy contracts. With regards to short-term operational time frame hedging against price fluctuation, a comparative analysis between forward contracts, call options and incentivizing end-consumers for flexibility provision is performed in \cite{Zhou2017}. In the context of our research, hedging is defined to be similar to the concept of Forward Contracts \cite{Kirschen2004} where the consumer wants to ensure that in the day-ahead market, they do not pay a price higher than their maximum willingness to pay. To illustrate this, assume that the consumers maximum willingness to pay is \euro70/MWh. For determining the price of electricity, we use the concept of Locational Marginal Price (LMP) which comprises of the marginal cost of generation and cost of network congestion \cite{Rivier1993}. If the LMP is \euro80/MWh, the consumer then to ensure that they pay at their maximum willingness to pay would need to offer energy flexibility through demand-side management. The question that we answer in this paper is: \textit{How can demand-side flexibility be used to hedge against price volatility in a distribution grid with distributed generation?} To address this question, we use a concept presented in \cite{Baker2016} to constrain the marginal prices to the maximum willingness of a consumer to pay for electricity. \par



The paper provides novel contributions to the scientific community, by investigating the possibility of using flexibility through demand-side management to hedge against electric price volatility. While there exists limited research on the applicability of demand-side flexibility as an alternative hedging mechanism, an approach that is able to quantify the energy flexibility required while accounting for operational constraints such as grid congestion are lacking.    


The hedging mechanism proposed in this paper can also be extended to the transmission grid. However, large industrial consumers in this grid, generally have long-term bilateral contracts with their energy provider and are hence protected against price volatility. Our focus in this paper is on the distribution grid and the entities connected to it such as small-scale industrial, commercial or residential consumers. The hedging mechanism presented in this paper could provide beneficial to the distribution grid based market entities, especially under the assumption of future institutional arrangements where dynamic pricing is adopted in distribution grid markets to enable more demand response and consumer participation. \par

The second contribution of our paper is an Organizational Structure for the coordination between the Distribution System Operator (DSO) and aggregator for demand-side flexibility management at the distribution grid. Lack of coordination between DSO and aggregators can have negative impacts on the operation of the distribution grid. As an example, research conducted in \cite{Verzijlbergh2014} illustrates that profit-seeking aggregators that try to benefit from price fluctuations can inadvertently cause congestion in the distribution grid. Our proposed approach is able to provide insights on how energy flexibility management through DSO-aggregator coordination can hedge against price volatility and price spikes in electricity markets.

The structure of this paper is as follows: first we present the problem formulation using duality theory to constrain dual variables and its impact on the primal problem. We then explain its application to the optimal power flow problem and the institutional arrangement between the DSO and aggregator for flexibility management in distribution grids. An illustrative case study is then presented from which insights on the business opportunity for aggregators are drawn. Finally we summarize our findings and give recommendations for future research.

\section{Power Flow Formulation for Hedging Against Price Volatility using Demand-Side Flexibility}\label{sec3}

In this Section, we illustrate the impact of constraining dual variables on the primal problem in the setting of an optimal power flow formulation. For doing so, we will build on the formulation that was previously introduced in \cite{Baker2016}.  We assume that the consumers, in order to hedge against the price volatility can use demand-side flexible resources. The formulation for determining the amount of demand-side flexibility required for constraining prices is presented in the following sub-sections. \par 

\subsection{Economic Dispatch Formulation}\label{subsec31}

In Economic Dispatch, all generators and loads are considered to be connected to the same bus, thereby ignoring physical network constraints. We assume the consumer derives utility of \euro b/MWh from load consumption and pays \euro a/MWh for it. As an entity, the customer would prefer to maximize its utility function through flexible consumption of electricity, and the upper and lower bounds for this are defined by \(\overline{P_L}\) and \(\underline{P_L}\) respectively. Neglecting network losses, this formulation is expressed as follows:  
\begin{equation}\label{eq:SWEDGeneral}
\begin{aligned}
& \underset{P_G, P_L}{\text{maximize}}
& & bP_L - aP_G \\
& \text{subject to}
& & \underline{P_L} \leq P_L \leq \overline{P_L} \qquad:\mu_1, \mu_2\\
& & & P_G - P_L = 0 \: \enspace \qquad :\lambda\\
& & & P_G, P_L \geq 0
\end{aligned}
\end{equation}
Dual variables \(\mu_1\), \(\mu_2\) and \(\lambda\) associated with the load consumption limits and power balance are introduced. As the primal problem is linear in nature, the problem has convex characteristics, and Slater's condition will hold \cite{Boyd2004}. This implies that the optimal value of the dual problem is equal to that of the primal problem. We now introduce a constraint on the LMP such that it has to be less than or equal to the maximum willingness of a consumer to pay for electricity, \(\pi^{des}\). The modified dual problem is expressed in \eqref{eq:DualSWEDGeneral}.    

\begin{equation}\label{eq:DualSWEDGeneral}
\begin{aligned}
& \underset{\mu_1, \mu_2, \lambda}{\text{minimize}}
& & -\underline{P_L}\mu_1 + \overline{P_L}\mu_2 \\
& \text{subject to}
& & -\mu_1 + \mu_2 - \lambda \geq b \\
& & & \lambda \leq \pi^{des} \\
& & & \mu_1, \mu_2, \pi^{des} \geq 0
\end{aligned}
\end{equation}

 Adding a constraint on the LMP (\(\lambda\)) in the dual formulation, introduces a new variable in the primal problem. This new variable is defined as \(P^{flexreq}\), and is the amount of demand-side flexibility required to constrain the LMP to the consumer's maximum willingness to pay, i.e. \(\lambda \leq \pi^{des}\). The modified primal problem is:
\begin{equation}\label{eq:ModifiedSWEDGeneral}
\begin{aligned}
& \underset{P_G, P_L, P^{flexreq}}{\text{maximize}}
& & bP_L - aP_G - \pi^{des}P^{flexreq} \\
& \text{subject to}
& & P_G - P_L + P^{flexreq} = 0 \\
& & & \underline{P_L} \leq P_L \leq \overline{P_L} \\
& & & P_G, P_L, P^{flexreq} \geq 0
\end{aligned}
\end{equation}

\subsection{Economic Dispatch formulation with Network Flows}\label{subsec32}

\begin{figure}[htp]
\centering
\vspace{-0.5em}
\includegraphics[width=5cm]{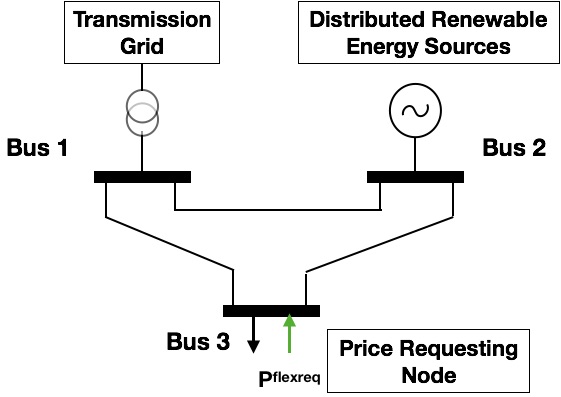}
\vspace{-.2em}
\caption{Simple grid for illustrating proposed hedging strategy}
\label{fig:ps2}
\vspace{-0.5em}
\end{figure}

The physical line capacity of the grid may further influence price volatility.
The physical structure of our distribution grid is shown in Figure \ref{fig:ps2}. This three bus system is taken for simplicity and the formulation can be extended for larger grids. 

In Figure \ref{fig:ps2}, at Bus 1 the distribution grid is connected to the transmission grid and it is considered to be the slack bus. Distributed energy resources are integrated to the distribution grid at Bus 2. The consumer who requests a price is referred to as a \textit{Price Requesting Load} and is connected at Bus 3, which is the price constrained bus. For modeling the cost of transmission grid we use the wholesale market prices, expressed as \(a^{trans}\) for the marginal cost and constant cost \(c^{trans}\) values. Similarly, the cost of distributed energy resources is expressed using marginal and constant costs as \(a^{dist}\) and \(c^{dist}\) respectively. The utility function of the price requesting load is comprised of marginal utility \(b^{load}\) and constant utility \(c^{load}\) values. The utility maximization formulation over 24 hours is as follows: 


\begin{subequations} \label{eq:DCOPFTransDist}
\begin{equation}
    \begin{split}
       \underset{P_G, P_L, P^{flexreq}}{\text{maximize}}
\sum^{24}_{l=1}b^{load}_l(P^{load})_l + c^{load}_l - a^{trans}_lP^{trans}_l \\- c^{trans}_l - a^{dist}_lP^{dist}_l - c^{dist}_l - \pi^{des}_lP^{flexreq}_l 
    \end{split}\tag{\ref{eq:DCOPFTransDist}}
\end{equation}
subject to:
\begin{equation}
(P_{G_{i_l}}) - (P_{L_{i_l}}) = \sum_{j \in \Omega_i}\frac{\theta_{i_l} - \theta_{j_l}}{X_{ij}}, \quad \forall i \in N \backslash K\label{eq:DCOPFTransDistA} 
\end{equation}
\begin{equation}
  (P_{G_{k_l}}) - (P_{L_{k_l}}) + (P^{flexreq}_{k_l}) = \sum_{j \in \Omega_k}\frac{\theta_{k_l} - \theta_{j_l}}{X_{kj}} \enspace \forall k \in K \label{eq:DCOPFTransDistB}
\end{equation}
\begin{equation}
    \theta_1 = 0 \label{eq:DCOPFTransDistC}
\end{equation}
\begin{equation}
    \underline{P_{L_i}} \leq (P_{L_{i_l}}) \leq \overline{P_{L_i}} \quad \forall i \in N \label{eq:DCOPFTransDistD}
\end{equation}
\begin{equation}
    -\overline{P_{ij}} \leq \frac{\theta_{i_l} - \theta_{j_l}}{X_{ij}} \leq \overline{P_{ij}} \quad \forall (i,j) \in \Omega_{ij} \label{eq:DCOPFTransDistE}
\end{equation}
\begin{equation}
    0 \leq (P_{G_{i_l}}) \leq \overline{P_{G_i}} \label{eq:DCOPFTransDistF} \quad \forall i \in N
\end{equation}
\begin{equation}
    P^{flexreq}_{k_l} \geq 0 \quad \forall k \in K \label{eq:DCOPFTransDistG}
\end{equation}
\end{subequations}
Sets N, K, \(\Omega_i\), \(\Omega_k\) and \(\Omega_{ij}\) represent the set of all nodes in the system, the price constrained nodes, buses connected to Bus i and Bus k, and the network lines respectively. Equation \eqref{eq:DCOPFTransDistA} represents the general power balance across nodes, while Equation \eqref{eq:DCOPFTransDistB} expresses the power balance at the price requesting node where demand-side flexibility is provided. Through \eqref{eq:DCOPFTransDistC} we express Bus 1 as the angle reference bus. The minimum and maximum bounds on loads at each bus, the line flow limits, upper and lower bounds on generation and positivity constraint on required flexibility is stated through  Equations \eqref{eq:DCOPFTransDistD}-\eqref{eq:DCOPFTransDistG} respectively. 


\section{Organizational Structure}\label{section:Organizational Structure}

\begin{figure}[htp]
\vspace{-1em}
\centering
\includegraphics[width=6.5cm]{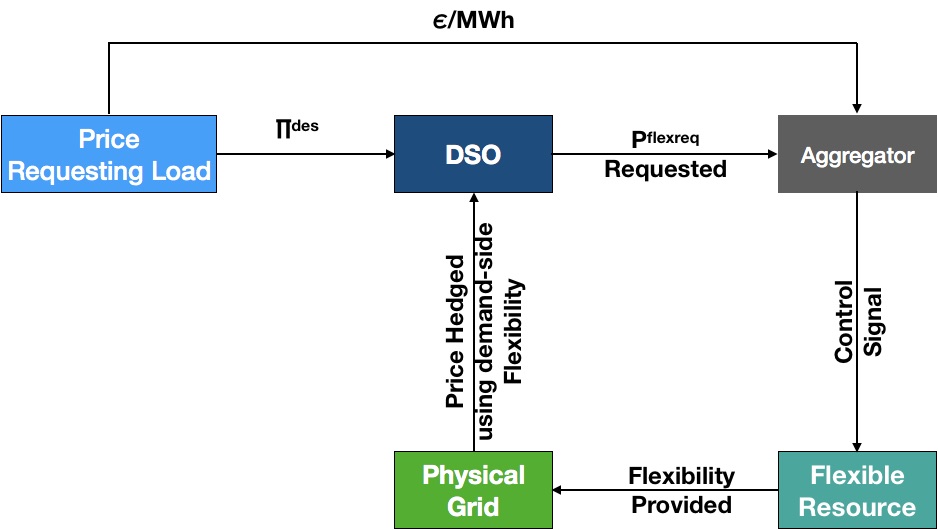}
\vspace{-.2mm}
\caption{Financial and Information Flow}
\label{fig:ps1}
\end{figure}

We propose an Organizational Structure, in association with the physical structure presented in Figure \ref{fig:ps2}, through which we focus on actor interaction and information flow. In Figure \ref{fig:ps1}, a \textit{Price Requesting Load} is situated at Bus 3 of Figure \ref{fig:ps2}, who communicates it's maximum willingness to pay for electricity, \(\pi^{des}\), to the DSO. It is expected that the DSO in a future electricity market will be a regulated market facilitator in the distribution grid \cite{VlerickBusinessSchool2015}. Thus, the DSO being the only entity having information about generation costs, the maximum willingness to pay of the price requesting load, as well as the physical structure of the network, is able to run the optimal power flow formulated in \eqref{eq:DCOPFTransDist}.
\vspace{-1em}

Executing the optimal power flow with \(\pi^{des}\) specified enables the DSO to determine the required demand-side flexibility for constraining the LMP. If the LMP is lower than \(\pi^{des}\), then no additional flexibility is required. However, if the LMP is greater than \(\pi^{des}\), then the DSO determines the amount of flexibility required (\(P^{flexreq} > 0\)). The DSO then, on behalf of the price requesting load, issues a request for demand-side flexibility to an aggregator. In response, the aggregator controls a demand-side flexible resource to provide the flexibility required to constrain the LMP to the \(\pi^{des}\) specified. The aggregator, is remunerated for its services by the price requesting load. In this paper, a simplistic remuneration scheme for the aggregator is proposed based on the difference in the LMP and specified \(\pi^{des}\), and the amount of flexibility that the aggregator provides to constrain the LMP. For a 24 hour period, the revenue at each hour is represented in \eqref{eq:AggregatorProfit}. 

\begin{equation}\label{eq:AggregatorProfit}
\begin{aligned}
R = \sum_{l = 1}^{24}(\lambda_l - \pi^{des}_l) \times P^{flexreq}_{k_l}
\end{aligned}
\end{equation}

For hour \(l\), \(\lambda_l\) represents the LMP without any demand-side flexibility and \(P^{flexreq}_{k_l}\) represents the amount of demand-side flexibility provided by the aggregator at the price requesting load bus. In practice, the aggregator's revenue is dependent on the amount of energy flexibility that the aggregator could provide at a given instance to constrain the prices. Furthermore, the aggregator's revenue should be less than the cost of dispatching energy from the next available and more expensive generator, otherwise the prices would not be constrained at the consumer specified value. In terms of risk, by investing in a demand-side flexible resource and providing flexibility, the Aggregator is taking over the risk of price volatility and price spikes from the price requesting load on to themselves. Hence, the aggregator's revenue is contingent on the amount of price volatility and price spikes in the electricity market.

\section{Simulation and Results}\label{section:SimAndResults}

For understanding how the presented optimal power flow formulation with constrained LMP determines the amount of flexibility required, we run an illustrative simulation using Equation \eqref{eq:DCOPFTransDist}. To perform our simulation, we use synthetic data based on the Amsterdam Power Exchange (APX) \cite{ENTSOE}, such that the variable \(\pi^{trans}\) represents a 24 hour time series data. The cost of distributed energy resources (\(a^{dist}\), \(c^{dist}\)) is assumed by scaling down the value of the wholesale market prices by \(30\%\) and then sampling 24 random points within the range. Similarly, the load utility values (\(b^{load}\), \(c^{load}\)) are set by sampling 24 random values within the values of \(a^{trans}\) and \(a^{dist}\). The network structure used in the following simulation corresponds to Figure \ref{fig:ps2} and has \(x_{12} = x_{13} = x_{23} = 0.1\)p.u. reactance with the maximum flows under infinite and finite line capacities as \(1.0\)MW and \(0.6\)MW respectively.      



\subsection{Price Hedging under Infinite Line Capacity}\label{section: PriceHedgingInfLine}

Power flows across the network under this simulation are not constrained and the marginal cost coefficient of the transmission system (\(a^{trans}\)) is higher than that of the distributed energy sources \(a^{dist}\). This implies that according to the merit order, the distributed energy resources will be dispatched first. However, owing to their intermittent nature, the upper limit on the amount of energy that can be drawn from them is variable, and may not be sufficient to supply the load by themselves. Hence, energy import from the transmission grid would be required to satisfy the residual demand, making the transmission grid based generation as the marginal generator. Due to infinite line capacity, the LMP will be the same across all the buses and equal to the cost of the marginal generator.\par

\begin{figure}[htp]
\vspace{-1em}
\centering
\includegraphics[width=10cm, scale =2]{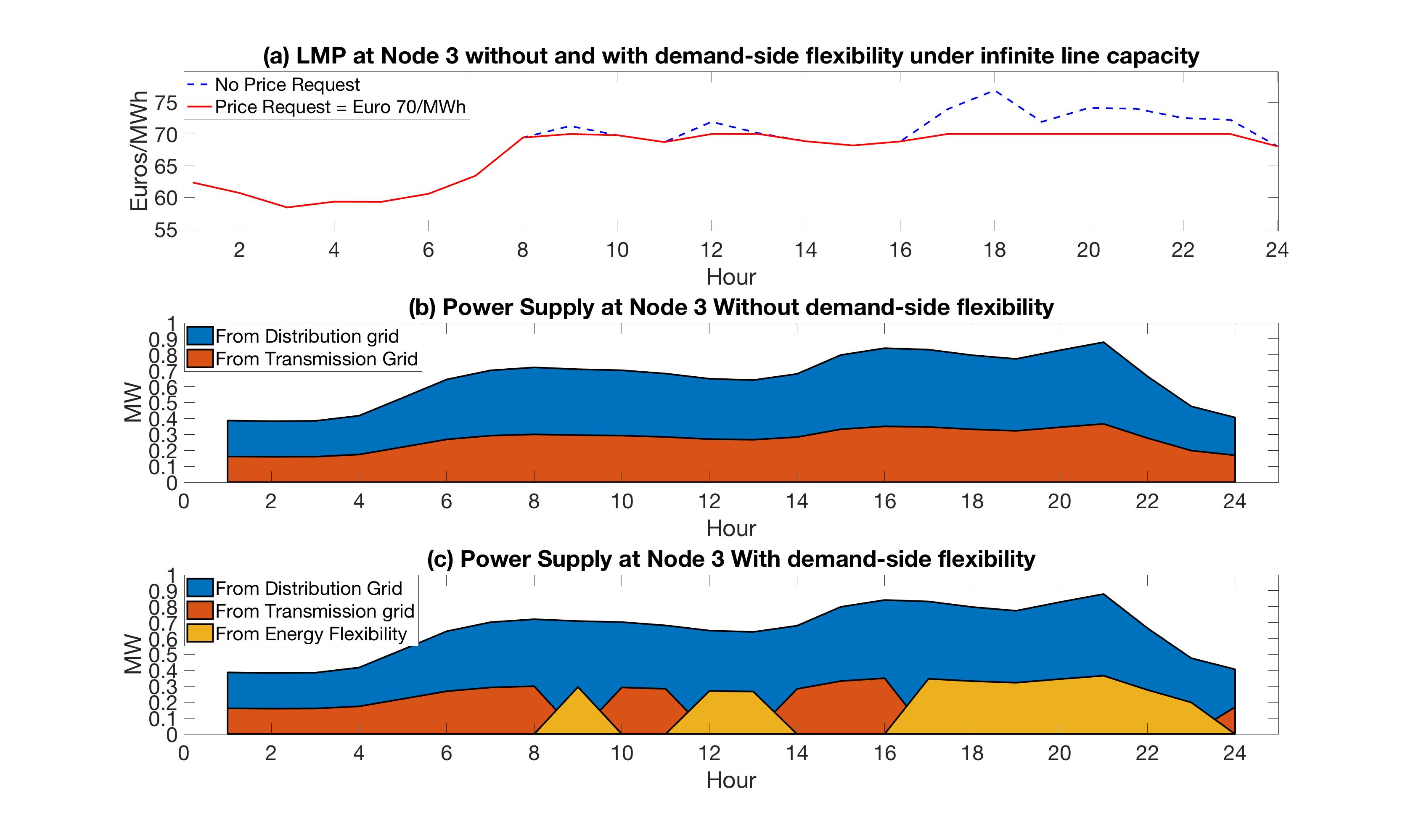}
\vspace{-2em}
\caption{Overview of constraining prices by specifying the Maximum Willingness to pay for Electricity under infinite line capacity: (a) LMP at Price Requesting Load with and without demand-side flexibility (b) Power supply without demand-side flexibility (c) Power supply with demand-side flexibility}
\label{fig:parameterizedBaseCaseLMP}
\vspace{-0.5em}
\end{figure}

The price requesting load at Bus 3, if it wants to hedge against the electricity price volatility, can specify its \mbox{maximum} willingness to pay for electricity directly to the DSO. DSO, having all the information it requires, can compute the amount of flexibility required for constraining the LMP. For our simulation purposes, we assume that the maximum willingness to pay for electricity is \(\pi^{des}\)=\euro70/MWh. 
It is observed from Figure \ref{fig:parameterizedBaseCaseLMP}, that when the LMP crosses the price threshold of \(\pi^{des}\)=\euro70/MWh at hour 9, the DSO would need to compute the flexibility required for constraining the LMP and communicate this information to the aggregator.

Figure \ref{fig:parameterizedBaseCaseLMP}(a) shows the LMP when the price is unconstrained and constrained. The quantification of demand-side flexibility is displayed through Figure \ref{fig:parameterizedBaseCaseLMP}(b) where we see the energy import from the expensive transmission grid based generation, and Figure \ref{fig:parameterizedBaseCaseLMP}(c) where demand-side flexibility can be provisioned to substitute this import thereby constraining prices. As an example, at hour 18, we observe that the proposed hedging mechanism can lead to a difference of \euro6.9/MWh. To achieve this hedging, the DSO computes and determines the flexibility required to be 0.33MW, which is then requested from the aggregator. The aggregator having the amount of flexibility required (\(P^{flexreq}\)) specified, is able to directly control a demand-side based flexible resource thereby provisioning flexibility. Hence, at hour 18, the price requesting load with \(\pi^{des}\) = \euro70/MWh is able to hedge \euro6.9/MWh by using 0.33MW of demand-side flexibility and the aggregator would be remunerated by \euro2.28 for that hour.

\subsection{Price Hedging under Finite Line Capacity}\label{section: PriceHedgingFiniteLine}

The main difference in this case, as compared to Section \ref{section: PriceHedgingInfLine}, is the impact of finite power flow capacity. Line congestion limit between Bus 2 and Bus 3 is set to 0.6MW. This implies that during periods of line congestion, even though the load at Bus 3 could be satisfied by the cheaper distributed energy resources, it would not be feasible due to the binding physical constraint. Furthermore, the LMP value now comprises of an additional network congestion cost. This manifests itself as price spikes that are observed at Bus 3. These price spikes are undesirable from a consumer's perspective. 


\begin{figure}[htp]
\vspace{-1em}
\centering
\includegraphics[width=10cm, scale =2]{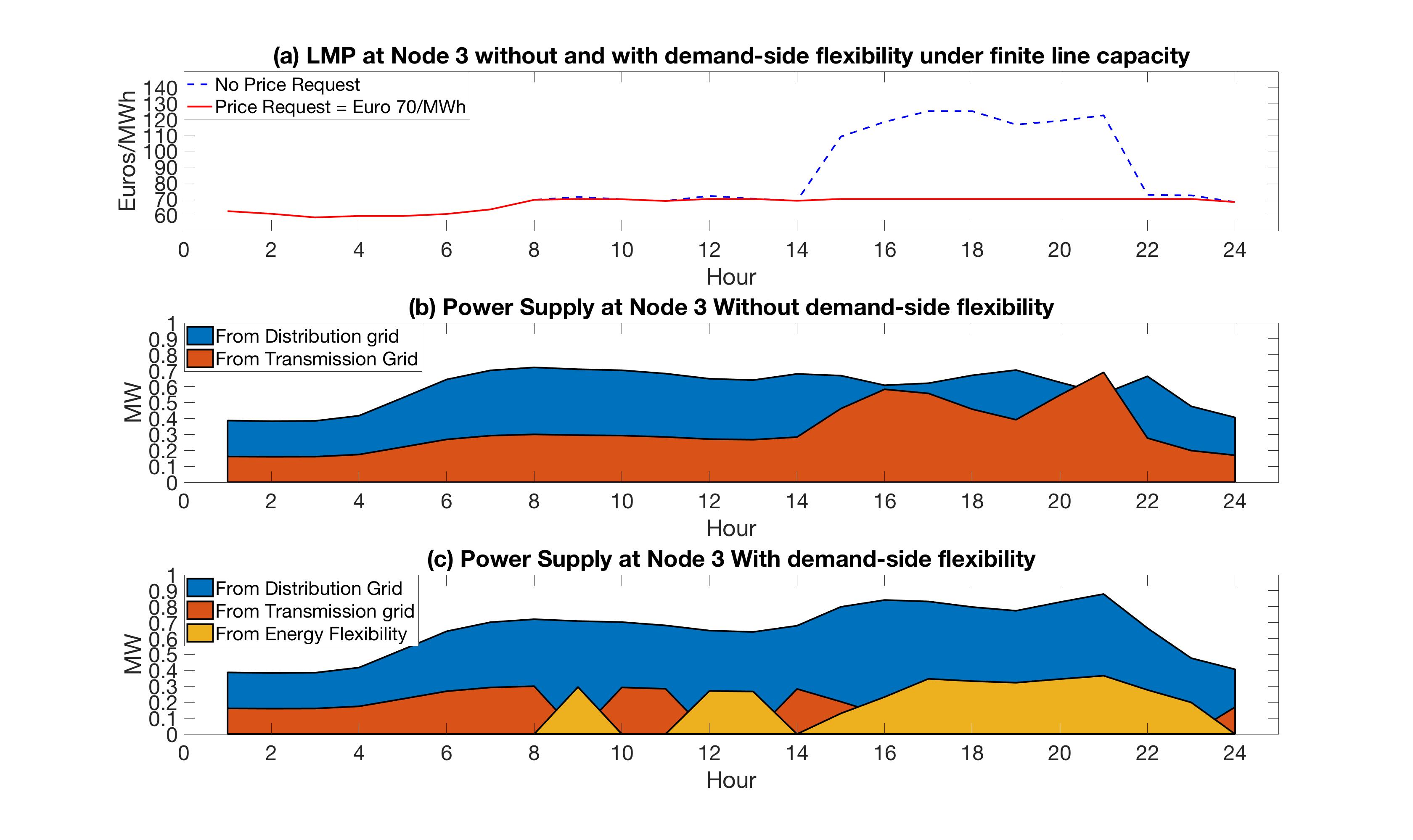}
\vspace{-2em}
\caption{Overview of constraining prices by specifying the Maximum Willingness to pay for Electricity under finite line capacity: (a) LMP at Price Requesting Load with and without demand-side flexibility (b) Power supply without demand-side flexibility (c) Power supply with demand-side flexibility}
\label{fig:priceVolFlexRequired}
\vspace{-0.5em}
\end{figure}

In Figure \ref{fig:priceVolFlexRequired}(a), price spikes are observed during the period of line congestion from hours 14-22. As observed in Figure \ref{fig:priceVolFlexRequired}(b), this period corresponds to that of increased energy import from the transmission grid which is more expensive. By specifying their maximum willingness to pay for electricity \(\pi^{des}\) to the DSO and eventually availing flexibility from the aggregator, the price requesting load is able to hedge against the price spikes. For constraining the prices, the amount of flexibility that is required is computed by the DSO and is depicted in Figure \ref{fig:priceVolFlexRequired}(c). At hour 17 for example, by provisioning 0.35MW of flexibility, the aggregator can cause a price difference of \euro55.14/MWh.The aggregator for its services would receive a remuneration of \euro19.3 for that hour, an amount which is significantly larger as compared to price hedging in the absence of price spikes. 

\subsection{Aggregator Revenue}\label{section: AggregatorRevenue}

In Section \ref{section:Organizational Structure}, we expressed the revenue that the aggregator could receive from a price requesting load for demand-side flexibility over a 24 hour period. Using our presented formulation, we determined the amount of flexibility that the aggregator would need to provide for hedging against price volatility and price spikes. Using Table \ref{tab:table1}, we calculate the sensitivity of aggregator's revenue to the maximum willingness to pay of a price requesting load.

\begin{table}[h!]
  \begin{center}
    \caption{Aggregator Revenue under different \(\pi^{des}\) values}
    \label{tab:table1}
    \begin{tabular}{l|S|r} 
      \textbf{\(\pi^{des}\)} & \textbf{Uncongested} & \textbf{Congested}\\
      \hline
      \euro68/MWh & \euro17.03 & \euro129.44\\
      \euro70/MWh & \euro9.22 & \euro106.97\\
      \euro72/MWh & \euro3.95 & \euro89.67\\
    \end{tabular}
  \end{center}
\end{table}

From Table \ref{tab:table1}, we see that as the difference between the LMP and \(\pi^{des}\) increases, the aggregator generates more revenue. Additionally, price spikes due to grid congestions can further increase money inflow for the aggregator. Given the results, it can be stated that there exists a business opportunity for aggregators.

\section{Conclusion}
In this paper, we presented a formulation to estimate the amount of demand-side flexibility that is required to hedge against electricity price volatility and price spikes caused by grid congestions. An organizational structure between DSO, aggregator and price requesting load for coordinating flexibility in distribution grids is proposed. By constraining the marginal price to the consumer's maximum willingness to pay for electricity, we use duality theory in quantifying the amount of flexibility required to perform hedging under infinite and finite line capacities. Accompanying simulation studies illustrate how line flows, generation profiles and flexible demand are affected by the proposed hedging mechanism. \par

These results can be extended in several new directions. Different contractual arrangements along with quantification of risk transfer between the DSO, aggregator and price requesting load for flexibility management will be explored. The results presented here do not account for inter-temporal constraints and losses in the power system and this will also be a focus of future research. Finally, a financial feasibility analysis of the proposed approach, in which we account for the cost of flexible resources, is required. 

\section*{Acknowledgment}
This work has received funding from the European Union's Horizon 2020 research and innovation programme under the Marie Sklowdowska-Curie grant agreement No. 675318 \mbox{(INCITE)}.



%

\bibliographystyle{ieeetr}
\bibliography{BibliographyUpdated}



\end{document}